\documentclass[11pt]{article}

\usepackage{amsmath}
\usepackage{graphicx}
\usepackage{epsfig}
\usepackage{amsfonts}
\usepackage{amssymb}
\usepackage{placeins}
\usepackage{cases}
\usepackage[margin=1in]{geometry}
\usepackage{authblk}
\usepackage[titletoc,toc,title]{appendix}
\usepackage{epstopdf}
\usepackage{xcolor}
\usepackage[final]{changes}
\usepackage{marvosym}
\usepackage{enumerate} 
\usepackage{bbm}

\usepackage{accents}

\title{Nonlinear Data-Driven Approximation of the Koopman Operator}  

\begin{document}
\author[1]{Dan Wilson \thanks{corresponding author:~dwilso81@utk.edu} }
\affil[1]{Department of Electrical Engineering and Computer Science, University of Tennessee, Knoxville, TN 37996, USA}
\maketitle

\begin{abstract}
Koopman analysis provides a general framework from which to analyze a nonlinear dynamical system in terms of a linear operator acting on an infinite-dimensional observable space.  This theoretical framework provides a rigorous underpinning for widely used dynamic mode decomposition algorithms.  While such methods have proven to be remarkably useful in the analysis of time-series data, the resulting linear models must generally be of high order to accurately approximate fundamentally nonlinear behaviors.  This issue poses an inherent risk of overfitting to training data thereby limiting predictive capabilities.  By contrast, this work explores strategies for nonlinear data-driven estimation of the action of the Koopman operator.   General strategies that yield nonlinear models are presented for systems both with and without control.   Subsequent projection of the resulting nonlinear equations onto a low-rank basis yields a low order representation for the underlying dynamical system.  In both computational and experimental examples considered in this work, linear estimators of the Koopman operator are generally only able to provide short-term predictions for the observable dynamics while comparable nonlinear estimators provide accurate predictions on substantially longer timescales and replicate infinite-time behaviors that linear predictors cannot.  
\end{abstract}

%

\section{Introduction}
Model identification is a necessary first step in the design, optimization, control, and estimation of complex dynamical systems.  When the mechanisms that underlie the dynamics are well understood, models can often be derived using first principles approaches and subsequently fit to available data.   However, in applications where an underlying system is too complicated to write down the underlying equations, data-driven model identification can be a powerful alternative \cite{kutz16}, \cite{brun19}.  Substantial progress has been made in recent years in the development of algorithms for inferring dynamical models strictly from time-series data.   Dynamic mode decomposition (DMD) \cite{kutz16}, \cite{schm10}, \cite{rowl09} is one such algorithm, with the ability to represent the evolution of snapshot data in terms of a collection of linear modes with associated eigenvalues that determine the growth/decay/oscillation rates.  This general framework has inspired numerous variations that can, for instance, incorporate the influence of an exogenous control input \cite{proc16}, account for noise and uncertainty \cite{hema17}, and continuously adjust when the system parameters are time-varying \cite{zhan19}.  

While DMD has been used in a wide variety of applications to explicate the underlying behavior of snapshot data in terms of eigenmode and eigenvalue pairs, without additional modifications it yields a linear estimator for the underlying dynamics.  Alternative approaches have been developed to identify fully nonlinear representations for the underlying equations from data \cite{mang19}, \cite{brun16b}, \cite{pant19}, \cite{rudy17}.  These methods typically consider a large nonlinear function library and subsequently use machine learning algorithms to choose a sparse subset that best matches the training data.  While such methods can be readily applied to identify sparse representations of highly nonlinear (and even chaotic systems), their efficacy is dependent on the choice an appropriate nonlinear function library.  Related machine learning approaches using neural networks have also achieved success for prediction in highly nonlinear dynamical systems \cite{path18}, \cite{vlac18}, \cite{rais18}.  

Many of the data-driven model identification strategies described above have a close connection to Koopman analysis \cite{budi12}, \cite{mezi13}, \cite{mezi19}.  Koopman-based approaches generally allow for the representation of a nonlinear dynamical system as a linear operator acting on an infinite-dimensional observable space.   Such approaches are distinct from standard linearization techniques that consider the dynamics in a close neighborhood of some nominal solution.  Rather, the goal of Koopman analysis is to identify a linear operator that can accurately capture fundamentally nonlinear dynamics -- the key challenge  is in the identification of a suitable finite basis to represent the action of the generally infinite dimensional Koopman operator.  The connection between DMD and spectral analysis of the Koopman operator is well established \cite{rowl09}, and in applications where high-dimensional data is readily available, the DMD algorithm can indeed be used to provide a finite-dimensional approximation of the Koopman operator.  Extensions of the DMD algorithm have illustrated that more accurate approximations of the Koopman operator can be obtained using DMD in conjunction with a set of lifting functions \cite{will15} and/or time-delayed embeddings of snapshot data \cite{arba17}.   Additional accuracy can also be obtained using an adaptive strategy that uses DMD to obtain a continuous family of linear models and actively chooses the one that provides the best representation at every instant \cite{wils22data}.    

How best to approximate the action of the Koopman operator from snapshot data remains an open question.  DMD separates data into snapshot pairs and subsequently finds a linear operator that provides a least squares fit for the mapping from one snapshot to the next.  This is currently the most widely used approach.    An obvious advantage of linear estimators of the Koopman operator is that they allow for subsequent analysis using a wide variety of linear techniques.  Nonetheless, such linear estimators are not always suitable for highly nonlinear systems since finite dimensional linear operators cannot be used, for instance, to replicate the infinite-time behavior of systems with multiple hyperbolic fixed points or systems with stable limit cycles.  Further limitations of linear estimators can also be seen in \cite{page18} which established the difficulty of representing even the relatively simple Burgers' equation in terms of a linear operator owing to the existence highly degenerate Koopman eigenvalues.  Alternatively, nonlinear models obtained from data-driven techniques are often more difficult to analyze, but can admit lower dimensional realizations and can often provide accurate representations of chaotic behavior \cite{brun16b}, \cite{path18}, \cite{vlac18}.  Recent works have considered nonlinear estimation strategies.  For instance, \cite{piet19} approximates separate Koopman operators that result for different values of an applied control input and uses this information to  formulate a switching time optimization problem.  Related approaches consider bilinear approximations of the Koopman operator for control systems \cite{peit18}, \cite{peit20}, \cite{sura16}.

This work explores strategies for nonlinear data-driven estimation of the action of the Koopman operator.  General strategies that yield nonlinear models are presented for systems both with and without control.  In the various examples considered in this work, only short term predictions of the dynamical behavior can be obtained using linear estimators for the Koopman operator.  By contrast, nonlinear estimators are able to provide accurate long-term estimates for the dynamics of model observables and yield accurate information about limit cycling behaviors and basin of attraction estimates.  The organization of this paper is as follows:~Section \ref{koopbackground} provides necessary background on Koopman operator theory along with a brief description of associated data-driven model identification techniques including DMD \cite{kutz16}, extended DMD \cite{will15}, and Koopman model predictive control \cite{kord18}.  Section \ref{koopnonlin} proposes algorithms for obtaining a nonlinear approximation for the Koopman operator from snapshot data in both autonomous and controlled systems.  The proposed approach is related to the extended DMD algorithm in that it considers a dictionary of functions of the observables, however, instead of estimating the action of the Koopman operator on each of the elements of the dictionary, the explicit nonlinear dependence of the dictionary elements on the observables is retained.  A variety of examples are presented in Section \ref{exampsec}.  Here, linear estimators for the Koopman operator are generally able to provide short-term predictions for the dynamics of observables; comparable nonlinear estimators provide accurate predictions on substantially longer timescales and accurately identify infinite-time behaviors.  Concluding remarks and suggestions for extension are provided in Section \ref{concsec}.


\section{Background}   \label{koopbackground}

\subsection{Koopman Operator Theory}

Consider a discrete-time dynamical system
\begin{equation} \label{vecfield}
x^+ = F(x),
\end{equation}
where $x \in \mathbb{R}^n$ is the state and $F$ gives the potentially nonlinear dynamics of the mapping $x \mapsto x^+$.  The Koopman operator $K:\mathcal{F} \rightarrow \mathcal{F}$ acts on the vector space of observables so that 
\begin{equation}\label{koop2}
K \psi (x) \equiv   \psi(F(x)),
\end{equation}
for every $\psi: \mathbb{R}^n \rightarrow \mathbb{R}$ belonging to the space of observables $\mathcal{F}$.   This operator is  linear (owing to the linearity of the composition operator).   As such, it can be used to represent the dynamics associated with a fully nonlinear system.  Approaches that use Koopman analysis are distinct from standard linearization techniques that are only valid in a close neighborhood of some nominal solution.  Note that while the Koopman operator is linear, it is generally infinite-dimensional \cite{budi12}, \cite{mezi13}, \cite{mezi19}.  In practical applications, the critical challenge of Koopman analysis is in the identification of a finite-dimensional approximation of the Koopman operator.

\subsection{Finite Dimensional Approximation of the Koopman Operator}   \label{koopest}
Dynamic mode decomposition (DMD) \cite{kutz16}, \cite{schm10}, \cite{tu14} is one standard approach for identifying a finite dimensional approximation of the Koopman operator.  To summarize this algorithm, one can consider a series of data snapshots
\begin{equation}
s_i = (g(x_i),g(x_i^+)),
\end{equation}
for $i = 1,\dots,d$ where $g(x) \in \mathbb{R}^m = \begin{bmatrix} \psi_1(x) , \dots, \psi_m(x) \end{bmatrix}$ is a set of observables obtained from the data and  $x_i^+ = F(x_i)$.  The goal of DMD is to identify a linear dynamical system of the form
\begin{equation} \label{linmod}
g_i^+ = A g_i,
\end{equation}
where $g_i = g(x_i)$, $g_i^+ = g(x_i^+)$, and $A \in \mathbb{R}^{m \times m}$ maps the observables from one time step to the next.   Such an estimate can be found according to a least-squares optimization,
\begin{equation} \label{standarddmd}
A = X^+  X^\dagger,
\end{equation}
where $X \equiv [g_1 \dots g_d]$, $X^+ \equiv [g_1^+ \dots g_d^+]$, and $^\dagger$ denotes the pseudoinverse.  As a slight modification, instead of taking the pseudoinverse of $X$ as in Equation \eqref{standarddmd}, it is often desirable to obtain a lower rank representation by first taking the singular value decomposition of $X$ and truncating terms associated with low magnitude singular values \cite{proc16}, \cite{brun17}.  Notice that the DMD algorithm as described above does not require knowledge of the underlying state and as such, can be implemented in a purely data-driven setting.  DMD often struggles in applications where few observables are available, i.e.,~when $m$ is small.   In such cases, extended DMD (EDMD) can be used \cite{will15}, which considers a lifted observable space
 \begin{equation}  \label{hvec}
 h(x) = \begin{bmatrix}  g(x)   \\   f_{\rm lift}(g(x))   \end{bmatrix}   \in \mathbb{R}^{m+b},
 \end{equation}
where $f_{\rm lift}(g(x)) \in \mathbb{R}^b$ is a possibly nonlinear function of the observables called a `dictionary'.  As before, letting $h_i = h(x_i)$ and $h_i^+ = h(x_i^+)$ comprise snapshot pairs with 
\begin{align} \label{caph}
H & \equiv \begin{bmatrix} h_1 \dots  h_d \end{bmatrix}, \nonumber \\
H^+ & \equiv  \begin{bmatrix} h_1^+  \dots  h_d^+ \end{bmatrix},
\end{align}
an estimate for the Koopman operator using the lifted coordinates can be obtained according to $A_{\rm lift} = H^+ H^\dagger$.  The EDMD approach can provide more accurate estimates of the Koopman operator than the standard DMD approach.  Indeed, in some cases the estimated Koopman operator converges to the true Koopman operator in the limit as both the lifted state and number of measurements approach infinity \cite{klus16}, \cite{kord18b}.  Possible choices of lifted coordinates include polynomials, radial basis functions, and Fourier modes \cite{will15}.   Additionally, delay embeddings of time series measurements of observables \cite{brun17}, \cite{arba17} have also yielded useful results in a variety of applications.

\subsection{Koopman-Based Model Identification With Control}  \label{koopcontrol}
Koopman-based approaches can readily be generalized to actuated systems \cite{kord18, proc18, will16}.  Following the approach suggested in \cite{kord18}, consider a controlled dynamical system
\begin{equation} \label{nlinvec}
x^+ = F(x,u),
\end{equation}
with output also given by Equation \eqref{koop2}.  The above equation is identical to \eqref{vecfield} with the incorporation of a control input $u \in \mathbb{R}^q \subset \mathcal{U}$.  Following the approach from \cite{kord18}, one can define an extended state space that is the product of the original state space $\mathbb{R}^n$ and the space of all input sequences $l(\mathcal{U}) = \{  (u_i)_{i=0}^\infty |  u_i \in \mathcal{U}  \}$.  Defining an observable $\phi:\mathbb{R}^n \times l(\mathcal{U}) \rightarrow \mathbb{R}$ belonging to a space of observables $\mathcal{H}$, the nonautonomous Koopman operator $K:\mathcal{H} \rightarrow \mathcal{H}$ can be defined according to
\begin{equation}
K \phi (x,(u_i)_{i=0}^\infty) =   \phi(F(x,u_0),(u_i)_{i=1}^\infty).
\end{equation}
Leveraging the EDMD algorithm, an estimate for the nonautonomous Koopman operator can be obtained by defining a vector of lifted coordinates
\begin{equation}
p(x_i) = \begin{bmatrix} g(x_i) \\  f_{\rm lift}(g(x_i))  \\ u_i \end{bmatrix},
\end{equation}
and determining an estimate for the linear dynamical system $p(x_i^+) = A_c p(x_i)$, where $A_c \in \mathbb{R}^{(m+b+q) \times (m+b+q)}$.  As noted in \cite{kord18}, one is generally not interested in predicting the last $q$ components of $p(x_i^+)$, i.e.,~those associated with the control input.  As such, the estimation of the final $q$ rows of $A_c$ can be neglected.  Letting $\bar{A}$ correspond the first $m+b$ rows of $A_c$.  Partitioning $\bar{A} = \begin{bmatrix} A & B \end{bmatrix}$ with $A \in \mathbb{R}^{(m+b) \times (m+b)}$ and $B \in \mathbb{R}^{(m+b) \times q}$, a linear, finite dimensional approximation of the Koopman operator can be obtained using a series of snapshot triples
\begin{equation} \label{triplefn}
w_i = (h_i, h_i^+,u_i),
\end{equation}
for $i = 1, \dots, d$.  Recall that $h_i$ and $h_i^+$ were defined below Equation \eqref{hvec}.  Once again, defining $H$ and $H^+$ as in \eqref{caph} and letting $\Upsilon = \begin{bmatrix} u_1  \dots  u_d \end{bmatrix}$, an estimate for $\bar{A}$ can be obtained according to 
\begin{equation} \label{dmdcfit}
\bar{A} = \begin{bmatrix} A & B \end{bmatrix} = H^+ \begin{bmatrix}  H  \\  \Upsilon \end{bmatrix}^ \dagger,
\end{equation}
ultimately yielding the state space representation 
\begin{equation}  \label{statecontrol}
h_i^+ = A h_i + B u_i.
\end{equation}
Using the above equation, the evolution of the observables can be recovered from the first $m$ entries of $h(x)$.

\section{Nonlinear Approximations of the Koopman Operator} \label{koopnonlin}

\subsection{Nonlinear Predictors For Autonomous Systems} \label{autsys}

The estimation strategies summarized in Sections \ref{koopest} and \ref{koopcontrol} yield linear models, for instance, of the form \eqref{linmod} and \eqref{statecontrol}.  The strategy detailed below allows for additional nonlinear terms in the prediction of the dynamics.   To begin, consider an unperturbed, discrete time dynamical system of the form \eqref{vecfield} with observables $g(x) \in \mathbb{R}^m$.   Leveraging the delayed embedding approaches considered  in \cite{brun17} and \cite{arba17}, one can define a lifted state
\begin{equation} \label{liftstate}
\gamma_i = \begin{bmatrix}   h(x_i) \\ h(x_{i-1}) \\ \vdots \\ h(x_{i-z})   \end{bmatrix},
\end{equation}
where $z \in \mathbb{N}$ determines the length of the delayed embedding and $h(x)$ was defined in Equation \eqref{hvec}.   Here, $\gamma_i \in \mathbb{R}^{M}$ with $M = (z+1)(m+b)$.  Next, a secondary lifting is defined
\begin{equation}
\sigma_i = \begin{bmatrix}  \gamma_i \\ f_n(\gamma_i)   \end{bmatrix},
\end{equation}
where $f_n(\gamma_i) \in \mathbb{R}^L$ is an additional, generally nonlinear function of the lifted state $\gamma_i$.  The term $f_n$ represents an additional user specified lifting of the data.  For example, these terms can be comprised of polynomials, radial basis functions, and Fourier modes \cite{will15}.  \added{Letting $\sigma_i$ and $\sigma_i^+$ be the lifted coordinates on successive iterations}, a direct implementation of the EDMD algorithm detailed in Section \ref{koopest} would seek a matrix $A$ that solves
\begin{equation} \label{dmdminimization}
\min_A  \left[ \sum_{i=1}^d ||  \sigma_i^+  -   A \sigma_i    ||_F  \right],
\end{equation}
for a collection of data $(\sigma_i,\sigma_i^+)$ for $i = 1,\dots,d$ where $||\cdot||_F$ denotes the Frobenius norm.  Alternatively, one can instead neglect the prediction of the final $L$ states because they are direct functions of $\gamma_i$.  In this case, defining the matrix $\hat{A}$ to be the first $M$ rows of $A$ and letting  $\hat{A} = \begin{bmatrix} A_n & C_n  \end{bmatrix}$ where $A_n \in \mathbb{R}^{M \times M}$ and $C_n \in \mathbb{R}^{M \times L}$, the minimization problem becomes
\begin{equation} \label{optproblem}
\min_{A_n,C_n} \left[ \sum_{i=1}^d ||  \gamma_i^+  -   A_n \gamma_i -  C_n f_n(\gamma_i)    ||_F  \right].
\end{equation}
This minimization problem can be solved by computing
\begin{equation} \label{minsol}
\hat{A} = \begin{bmatrix} A_n & C_n  \end{bmatrix} = \Gamma^+ \begin{bmatrix} \Gamma \\ F_n  \end{bmatrix}^\dagger,
\end{equation}
where $\Gamma \equiv \begin{bmatrix} \gamma_1 \dots \gamma_d  \end{bmatrix}$,  $\Gamma^+ \equiv \begin{bmatrix} \gamma_1^+ \dots \gamma_d^+  \end{bmatrix}$, and $F_n = \begin{bmatrix} f_n(\gamma_1) \dots  f_n(\gamma_d)  \end{bmatrix}$.  The resulting model takes the form
\begin{equation} \label{nonlinpredict}
\gamma_i^+ = A_n \gamma_i + C_n f_n(\gamma_i).
\end{equation}
Lower rank approximations of $A_n$ and $C_n$ may be desirable in order to avoid overfitting to the measured data.  In this instance, one can consider the singular value decomposition 
\begin{equation} \label{svdmtx}
\begin{bmatrix} \Gamma \\ F_n  \end{bmatrix} = U \Sigma V^T,
\end{equation}
where $U \in \mathbb{R}^{(M+L) \times (M+L)}$, $\Sigma \in \mathbb{R}^{(M+L) \times d}$, and $V \in \mathbb{R}^{d \times d}$, and $^T$ denotes the matrix transpose.  Note that $U$ and $V$ are real because $\Gamma$ and $F_n$ are real.  A rank $r$ approximation of \eqref{svdmtx} can be obtained by taking letting $\tilde{U}$ and $\tilde{V}$ represent the first $r$ columns of $U$ and $V$, respectively, and letting $\tilde{\Sigma}$ be a square matrix containing the first $r$ singular values from $\Sigma$ so that 
\begin{equation}
\begin{bmatrix} \Gamma \\ F_n  \end{bmatrix} \approx \tilde{U} \tilde{\Sigma} \tilde{V}^T.
\end{equation}
With this representation, one can obtain the lower rank approximation of the solution of the optimization problem from \eqref{optproblem}
\begin{equation}
\begin{bmatrix} A_n & C_n  \end{bmatrix} \approx \Gamma^+ \tilde{V} \tilde{\Sigma}^{-1} \tilde{U}^T,
\end{equation}
where $^{-1}$ denotes the matrix inverse.

In contrast to the standard EDMD algorithm, the predictor \eqref{nonlinpredict} is nonlinear.  Nonetheless, as illustrated in the examples presented in Section \ref{exampsec}, the added nonlinearity can accommodate behaviors that linear predictors cannot.   

\subsection{Reduced Order Representations Using Nonlinear Predictors}  \label{redsec}

The model identification strategy proposed in Section \ref{autsys} incorporates a lifting of the observables in conjunction with a delayed embedding of the lifted coordinates.  As such, the resulting nonlinear model may be high dimensional making analysis and control difficult. Because the proposed strategy yields a nonlinear predictor for the dynamics of the observables of \eqref{vecfield} (as opposed to a linear predictor obtained from the EDMD algorithm), it is generally useful to identify a reduced order representation of the dynamics.  This task can be accomplished by applying proper orthogonal decomposition (POD) \cite{holm96}, \cite{rowl17} to $\Gamma$ to identify a representative set of modes from the data.  Here, POD modes are found according to the eigenvectors of $\Gamma \Gamma^T$ and sorted according to the magnitude of the associated eigenvalues.  Keeping the first $\rho$ POD modes and truncating the rest (i.e.,~associated with the smallest eigenvalues) yields an orthogonal basis of POD modes $\Phi = \begin{bmatrix} \mu_1 & \dots &\mu_\rho \end{bmatrix} \in \mathbb{R}^{M \times \rho}$ for which
\begin{equation} \label{gammaapprox}
\gamma_i \approx \sum_{k = 1}^\rho \mu_k \omega_{k,i},
\end{equation}
where $\omega_{k,i}$ is a coefficient that can be obtained according to $\omega_{k,i} = \mu_k^T \gamma_i$.  Substituting \eqref{gammaapprox} into \eqref{nonlinpredict} yields
\begin{equation}
\sum_{k = 1}^\rho \mu_k \omega^+_{k,i} \approx A_n \sum_{k = 1}^\rho \mu_k \omega_{k,i} + C_n f_ n \left( \sum_{k = 1}^\rho \mu_k \omega_{k,i} \right).
\end{equation}
Multiplying the above equation by the left by $\Phi^T$ and rearranging (noting that the POD modes are orthogonal) yields
\begin{equation} \label{lowdim}
\Omega^+_i = \Phi^T A_n \Phi \Omega_i + \Phi^T C_n f_n(\Phi \Omega_i),
\end{equation}
where $\Omega_i = \begin{bmatrix}  \omega_{1,i} & \dots & \omega_{\rho,i} \end{bmatrix}^T$ and $\Omega^+ = \begin{bmatrix}  \omega_{1,i}^+ & \dots & \omega_{\rho,i}^+ \end{bmatrix}^T$.  Equation \eqref{lowdim} provides an order $\rho$ approximation for the dynamics of the nonlinear system  given by Equation \eqref{nonlinpredict}.  Conversion from the reduced order basis $\Omega_i$ back to the lifted state $\gamma_i$ can be accomplished using Equation \eqref{gammaapprox}.

\subsection{Nonlinear Predictors For Controlled Systems} \label{nlincont}
Control input can readily be incorporated into the proposed model identification strategy.  To do so, considering a general system of the form \eqref{nlinvec},  one can define the lifted state as 
\begin{equation} \label{gammacont}
\gamma_{c,i} = \begin{bmatrix} h(x_i) \\ \vdots \\ h(x_{i-z}) \\ u_{i-1} \\  \vdots \\ u_{i-z}  \end{bmatrix}.
\end{equation}
Here, $\gamma_{c,i} \in \mathbb{R}^{M_c}$ where $M_c = (z+1)(m+b) + zq$.  Compared with the lifted state defined in Equation \eqref{liftstate},  $\gamma_{c,i}$ also contains an embedding of the preceding $z$ control inputs as suggested in \cite{arba18b}.  This lifted state is then augmented with additional states to yield
\begin{equation} \label{augstate}
\sigma_{c,i} = \begin{bmatrix} \gamma_{c,i} \\ u_i \\ f_{c,n}(\gamma_{c,i})  \\   \end{bmatrix},
\end{equation}
where $f_{c,n}(\gamma_{c,i}) \in \mathbb{R}^L$ is a nonlinear function of $\gamma_{c,i}$.  \added{Let $\sigma_{c,i}$ and $\sigma_{c,i}^+$ be the lifted coordinates on successive iterations}.  Mirroring the argument from Section \ref{autsys} that starts with Equation \eqref{dmdminimization} and ends with Equation \eqref{nonlinpredict}, for a collection of snapshot pairs $(\sigma_{c,i},\sigma_{c,i}^+)$ for $i = 1,\dots,d$,  a direct implementation of the EDMD algorithm would seek a matrix $A$ that solves the minimization problem $\min_A  \left[ \sum_{i=1}^d ||  \sigma_{c,i}^+  -   A \sigma_{c,i}    ||_F  \right]$.  However, prediction of the final $L+q$ states can be omitted because prediction of the control sequence is not of interest and $f_{c,n}(\gamma_{c,i})$ is an explicit function of $\gamma_{c,i}$.  Subsequently defining the matrix $\hat{A}$ to be the \added{first $M_c$ rows} of $A$ and letting $\hat{A} = \begin{bmatrix} A_c & B_c & C_c \end{bmatrix}$ the minimization problem becomes 
\begin{equation} \label{optcontrol}
\min_{\added{A_c,B_c,C_c}} \left[ \sum_{i=1}^d ||  \gamma_{c,i}^+  -   A_c \gamma_{c,i} -  B_c u_i  -  C_c f_{c,n}(\gamma_{c,i})    ||_F  \right],
\end{equation}
which can be solved by computing
\begin{equation}  \label{ahatest}
\hat{A} = \begin{bmatrix}  A_c & B_c & C_c  \end{bmatrix} = \Gamma_c^+ \begin{bmatrix}  \Gamma_c \\  U \\ F_{c,n} \end{bmatrix}^\dagger,
\end{equation}
where $\Gamma_c \equiv \begin{bmatrix} \gamma_{c,1} \dots \gamma_{c,d}  \end{bmatrix}$,  $\Gamma_c^+ \equiv \begin{bmatrix} \gamma_{c,1}^+ \dots \gamma_{c,d}^+  \end{bmatrix}$, $U = \begin{bmatrix} u_1 \dots  u_n  \end{bmatrix}$, and $F_{c,n} = \begin{bmatrix} f_{c,n}(\gamma_{c,1}) \dots  f_{c,n}(\gamma_{c,d})  \end{bmatrix}$.  The resulting model takes the form
\begin{equation} \label{contsys}
\gamma_{c,i}^+ = A_c \gamma_{c,i} + B_c u_i + C_c f_{c,n}(\gamma_{c,i}).
\end{equation}
As with the autonomous system of the form \eqref{nonlinpredict}, a lower rank approximation of the matrices $A_c$, $B_c$ and $C_c$ can be obtained using a truncated singular value decomposition.  Likewise, a reduced order model similar to Equation \eqref{lowdim} can be obtained by projecting \eqref{contsys} onto a reduced order basis of POD modes obtained from the data contained in $\Gamma_c$.

\section{Examples With Comparisons to Other Koopman-Based Approaches}  \label{exampsec}

\subsection{Forced Duffing Equation}
\begin{figure}[htb]
\begin{center}
\includegraphics[height=2.6in]{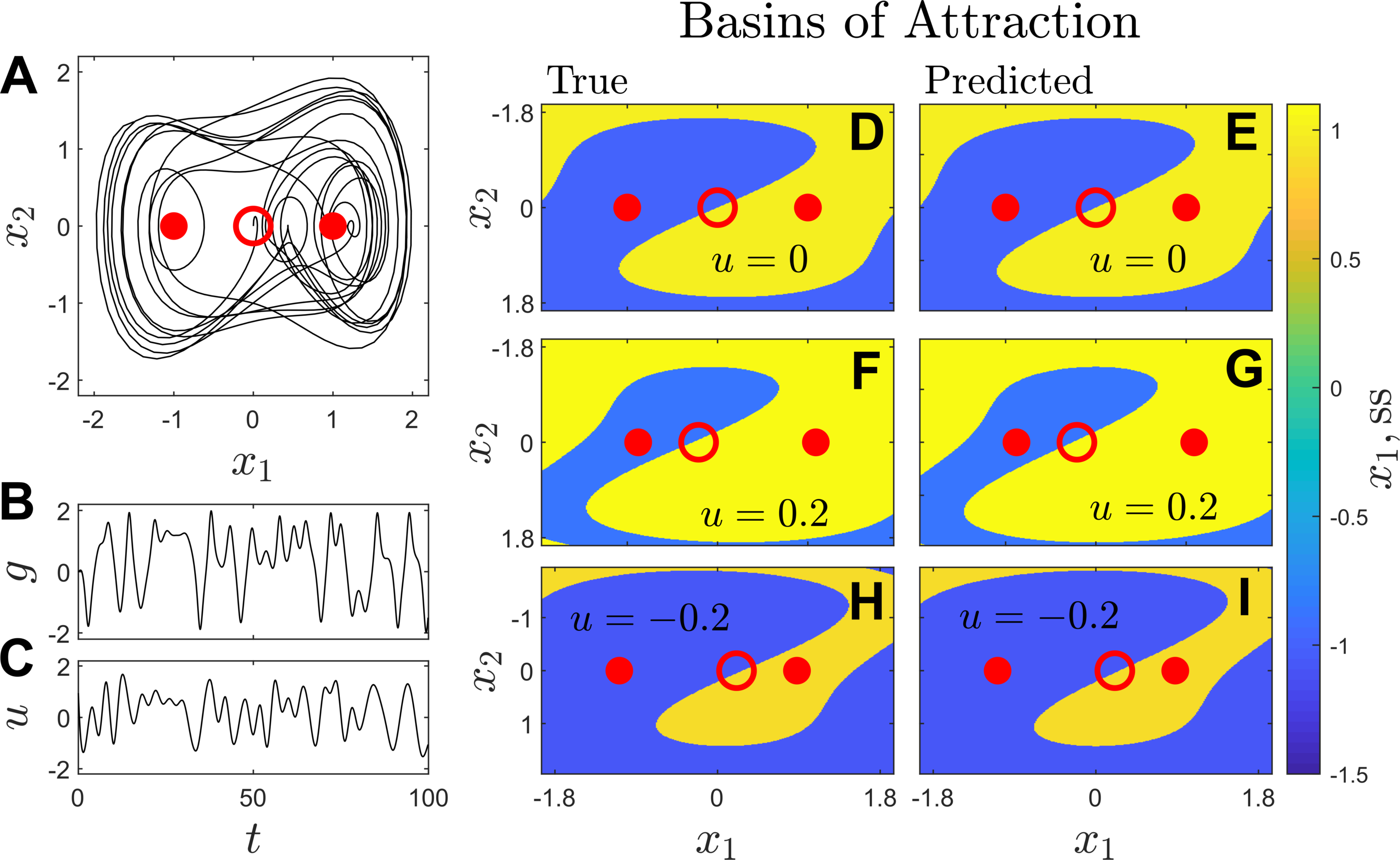}
\end{center}
\caption{Data-driven model identification and subsequent basin of attraction estimates for the forced Duffing Equation \eqref{duffingeq}.  Panel A shows a representative sample of the forced behavior of this model in response to the input from panel C.  Solid (resp.,~open) dots indicate stable (resp.,~unstable) fixed points.  Panels B gives the corresponding observable used for model identification.  Note that the state $x_2$ is not directly measured.  The model identification strategy from Section \ref{nlincont} is performed and the resulting nonlinear models are used to obtain basin of attraction estimates for the fixed points that result for various constant values of $u$.  Predicted basins of attraction are shown in panels E, G, and I when taking $u = 0, 0.2$, and $-0.2$, respectively.  Panels D, F, and H, show comparisons obtained from direct simulations of \eqref{duffingeq}.  Stable and unstable fixed points are shown for reference as closed and open circles, respectively.  The colormap represents the value of $x_1$ in the limit as time approaches infinity for both the inferred and actual models.      }
\label{duffingresults}
\end{figure}

Consider the forced Duffing equation
\begin{align} \label{duffingeq}
\dot{x}_1 &= x_2, \nonumber \\
\dot{x}_2 &= u(t) - \delta x_2 - \alpha x_1 - \beta x_1^3,
\end{align}
with observable
\begin{equation}
g(x) =  x_1,
\end{equation}
taking $\alpha = 1$, $\beta = -1$, and $\delta = 0.5$.  Here $u(t)$ represents a general control input instead of the usual periodic driving force.  When $u(t) = 0$, Equation \eqref{duffingeq} has one unstable equilibrium at $x_1 = x_2 = 0$ and two stable equilibria at $x_2 = 0$ and $x_1 = \pm 1$.   Data is obtained for model identification taking $u(t)$ as follows:~random numbers between -1.5 and 1.5 are chosen from a uniform distribution with the value held constant over a 5 time unit interval.  The resulting curve is smoothed with a spline interpolation and used as the input in Equation \eqref{duffingeq}.  Simulation is performed for $t \in [0, 1000]$ and the resulting output is used to implement the model identification procedure detailed in Section \ref{nlincont} taking snapshots at time intervals $\Delta t = 0.1$, i.e.,~so that $x_i = \begin{bmatrix}  x_1( \Delta t (i-1)) &  x_2( \Delta t (i-1)) \end{bmatrix}^T$.  Panel A of Figure \ref{duffingresults} shows the state of the system over the first 100 time units of simulation.  Panels B and C show the corresponding observable and input, respectively, used for model identification.  To implement the model identification strategy, a delay embedding of size $z = 1$ is used taking $h(x_i) = g(x_i)$ so that $\gamma_{c,i} \in \mathbb{R}^3$ as defined in Equation \eqref{gammacont}.  The nonlinear lifting $f_{c,n}(\gamma_{c,i}) \in  \mathbb{R}^{12}$ is comprised of polynomial terms in $h(x_i)$ and $h(x_{i-1})$ up to degree 4 (e.g.,~$h(x_i)^2, h(x_i)^2 h(x_{i-1}), h(x_i)h(x_{i-1}^3))$.  The matrix $\hat{A}$ is estimated according to Equation \eqref{ahatest} which is comprised of the matrices $A_c \in \mathbb{R}^{3\times 3}$, $B_c \in \mathbb{R}^{3\times 1}$, and $C_c \in \mathbb{R}^{3\times 12}$ from Equation \eqref{contsys}.

The inferred model is used to provide basin of attraction estimates for stable fixed points of the Duffing equation for different constant values of $u$.    For the inferred model, initial conditions are taken to be $\gamma_{c,1} = \begin{bmatrix} x_1  &  x_1 - \Delta t x_2  & 0 \end{bmatrix}^T$ and the associated basin of attraction is assigned according to the resulting steady state value of $x_1$, i.e.,~$x_{1,ss} = \lim_{j \rightarrow \infty} \left( e_1^T \gamma_{c,j} \right)$ where $e_1 = \begin{bmatrix} 1 & 0 & 0 \end{bmatrix}$.   Results are shown in panels D-I of Figure \ref{duffingresults}; basin of attraction estimates between the true model \eqref{duffingeq} and the inferred model of the form \eqref{contsys} are nearly identical.

For this example, comparison with linear estimators such as EDMD as described in Section \ref{koopcontrol} is not considered.  It is well known that linear models cannot be used to accurately represent the infinite time behavior of systems with multiple hyperbolic fixed points (as is the case in Equation \eqref{duffingeq}) complicating their use for providing basin of attraction estimates.  EDMD was used in \cite{will15} to obtain basin of attraction estimates of the unforced Duffing equation by considering the resulting approximation of the nontrivial Koopman eigenmode associated with nondecaying solutions.  This method of analysis, however, would require data from trajectories with initial conditions uniformly distributed over a domain of interest and with a constant value of $u$.  By contrast, the approximated basin of  attraction estimates from Figure \ref{duffingresults} are obtained from snapshot triples using arbitrary inputs and also provide accurate basin of attraction estimates for arbitrary values of $u$.

As a final note regarding this example using the forced Duffing equation, the model equations \eqref{duffingeq} can be recast in discrete time by first letting $\dot{x}_2(t) \approx (x_1(t)  - x_1(t - \Delta t))/\Delta t$.  Subsequently taking a forward Euler time step yields
\begin{align} \label{fealt}
x_1(t + \Delta t) &= a_1 x_1(t) + a_2x_1(t - \Delta t), \nonumber \\
x_2(t + \Delta t) &= a_3 u(t)  +a_4 x_1(t) +  a_5 x_1^3(t) + a_6x_1(t - \Delta t),
\end{align}
where $a_1 = 1, a_2 = -1, a_3 = \Delta t, a_4 = - (\delta + \alpha \Delta t), a_5 = -\beta \Delta t,$ and $a_6 = \delta$.  When using the delay embedding and polynomial lifting strategy described above, the augmented state \eqref{augstate} contains all of the polynomial terms that comprise \eqref{fealt}.   As such, for $\Delta t$ small enough, this example could readily be handled by a sparse nonlinear model identification algorithm \cite{brun16c}, \cite{brun16b} that selects appropriate functions from a library and identifies the associated coefficients.  The examples to follow, however, do not admit simple, sparse representations for the dynamics of the observables.

\FloatBarrier

\subsection{Conductance-Based Neural Model}
Consider a conductance-based Wang-Buzsaki model neuron \cite{wang96} with an additional adaptation current \cite{erme98}
\begin{align} \label{wbmodel}
C \dot{V}  &= -g_{\rm Na} m_\infty^3 p (V -E_{\rm Na}) - g_{\rm K} n^4(V -E_K) - g_{\rm L}(V-E_{\rm L}) - i_w  + u(t) + i_b, \nonumber \\
\dot{p} &= \gamma \left[ \alpha_p(V)(1-p) - \beta_p(V)p   \right], \nonumber \\
\dot{n} &= \gamma \left[ \alpha_n(V)(1-n) - \beta_n(V)n \right], \nonumber \\
\dot{w} &= a(1.5/(1+\exp((b-V)/k))-w). \nonumber \\
\end{align}
Here, $V$ is represents the transmembrane voltage with $p$ and $n$ representing gating variables.  The adaptation current $i_w = g_w w (V-E_K)$ is mediated by the variable $w$ and  $i_b = 10 \mu {\rm A}/{\rm cm}^2$ is a constant baseline current.  The input $u(t)$ represents a transmembrane current.  The membrane capacitance, $C$, is taken to be $1 \mu  {\rm F}/{\rm cm}^2$.  Auxiliary equations governing ionic currents are:
\begin{align*}
m_\infty &= \alpha_m(V)/(\alpha_m(V) + \beta_m(V)),   \\
\beta_n(V) &= 0.125\exp(-(V+44)/80),  \\
\alpha_n(V) &= -0.01(V+34)/(\exp(-0.1(V+34))-1), \\
\beta_p(V) &= 1/(\exp(-0.1(V+28))+1), \\
\alpha_p(V) &= 0.07\exp(-(V+58)/20), \\
\beta_m(V) &= 4\exp(-(V+60)/18), \\
\alpha_m(V) &= -0.1(V+35)/(\exp(-0.1(V+35))-1). \\
\end{align*}
Reversal potentials and conductance are $E_{\rm Na} = 55 {\rm m}V, E_{\rm K} = -90{\rm m}V, E_{\rm L} = -65 {\rm  m}V, g_{\rm Na}= 35 {\rm mS}/{\rm cm}^2, g_{\rm K} = 9 {\rm mS}/{\rm cm}^2, g_{\rm L} = 0.1 {\rm mS}/{\rm cm}^2, g_w = 2 {\rm mS}/{\rm cm}^2$. Auxiliary parameters are $a = 0.02\;{\rm ms}^{-1}, b = -5 \; {\rm m}V, k = 0.5 {\rm m}V$, and $\gamma = 5$.  In the absence of input, the neural model \eqref{wbmodel} is in a tonically firing regime with a stable limit cycle having period 6.53 ms. The input $u(t)$ serves to modulate the firing rate of the action potentials.

For the conductance-based neural model from Equation \eqref{wbmodel}, the state is $x = \begin{bmatrix} V & p & n & w \end{bmatrix}^T$.  The observable for the spiking neural model is taken to be 
\begin{equation}
g(x) = \begin{bmatrix} V \\ h \end{bmatrix},
\end{equation}
i.e.,~it is assumed that the variables $V$ and $h$ can be measured directly but that the variables $n$ and $w$ are inaccessable.  The model identification strategy from Section \ref{nlincont} is implemented using 300 ms of simulated data taking a time step of $\Delta t = 0.025$ ms with an applied input $u(t) = 6 \sin (2 \pi t/200 + 0.0003 t^2)$.  A delayed embedding of size $z = 10$ is used taking $h(x_i) = g(x_i)$ so that $\gamma_{c,i} \in \mathbb{R}^{32}$ as defined in Equation \eqref{gammacont}.  The nonlinear lifting function $f_{c,n}(\gamma_{c,i}) =  f_1(f_2(h(x_i)))$.  Here $f_2(h(x_i)) \in \mathbb{R}^{10}$ with the $j^{\rm th}$ term given by $||g(x_i) - q_j||_2$,  where $q_j \in \mathbb{R}^2$ is the center of each radial basis function with the first element (associated with the transmembrane voltage) chosen randomly from a uniform distribution taking values between -300 and 200 and the second element (associated with the gating variable) chosen randomly from a uniform distribution taking values between 0 and 1, and $||\cdot||_2$ denotes the 2-norm.  The function $f_1(f_2(h(x_i))) \in \mathbb{R}^{990}$ provides a second nonlinear lifting and is comprised by taking polynomial combinations of the elements of $f_2(h(x_i))$ up to degree 4.  The matrix $\hat{A}$ is estimated according to Equation  \eqref{ahatest} using a truncated singular value decomposition of rank 80 to approximate the pseudoinverse.  This information is used to determine $A_c \in \mathbb{R}^{32 \times 32}$, $B_c \in \mathbb{R}^{32\times 1}$, and $C_c \in \mathbb{R}^{32\times 990}$ from Equation \eqref{contsys}.  As described in Section \ref{redsec}, a 20 dimensional model is obtained by projecting the inferred model equations onto a POD basis obtained from the eigenvectors of $\Gamma_c \Gamma_c^T$.

\begin{figure}[htb]
\begin{center}
\includegraphics[height=2in]{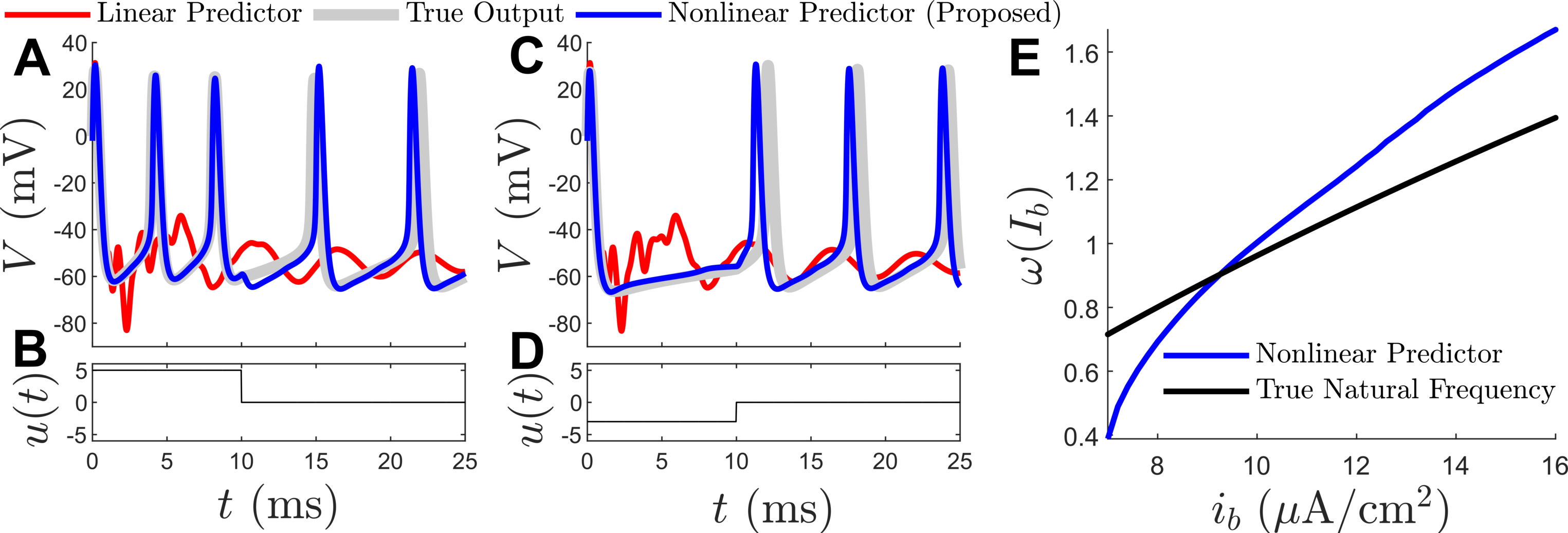}
\end{center}
\caption{Comparisons between the full order model, linear predictor, and proposed nonlinear predictor in response to various inputs.  Panel A (resp.,~C) shows the response to the pulse input in panel B (resp.,~D).  The proposed nonlinear predictor accurately reflects the increase (resp.,~decrease) in firing rate in response to positive (resp.,~negative) inputs as well as the subsequent rebound caused by the adaptation current. The linear predictor obtained using the strategy proposed in \cite{kord18} is unable to replicate the true model output.  Panel E highlights accurate predictions for the steady state firing rate predicted by the nonlinear inferred model for various baseline currents. Note that the inputs considered in the results shown above are substantially different than the oscillatory inputs used for model inference. }
\label{neuralresults}
\end{figure}

Simulations of the inferred model are compared to those of the simulations of the true model \eqref{wbmodel}.  Comparisons are also given when using the Koopman model predictive control strategy from \cite{kord18} which provides a least squares estimate for the update rule $a_i^+ = A a_i + B u_i$ where the lifted state space in this example is taken to be $a_i = \begin{bmatrix} \gamma_{c,i}^T & f_{c,n}(\gamma_{c,i})^T \end{bmatrix}^T$.  Results are shown in Figure \ref{neuralresults}.  Panel A shows the effect of a 10 ms duration positive pulse input shown in panel B.  Panel C shows the effect of a comparable negative pulse input from panel D.  In each case the proposed method (with the nonlinear predictor) provides a good approximation for the true model output while the linear predictor does not.  The model obtained from the nonlinear predictor yields stable oscillations in response to constant inputs -- such stable oscillations are not possible to obtain when considering linear predictors.  Pane E shows the predicted natural frequency for different baseline currents illustrating good agreement with the true model.  These results are particularly noteworthy considering that the model was trained using only oscillatory inputs and that this model was inferred without direct access to the auxiliary variables $n$ and $w$.

\subsection{Coupled Population of Neural Oscillators} \label{neurpopsec}
In the previous section, the dynamics of a single conductance-based neuron was considered.  Here, the behavior of a coupled population of identical, noisy thalamic neural oscillators taken from \cite{rubi04} will be considered:
\begin{align} \label{neurmod}
C \dot{V}_j &= f_V(V_j,p_j,r_j) + i_b + u(t) + \sqrt{2D} \eta_j + \frac{1}{N} \sum_{i = 1}^N  \sigma_c  (V_i - V_j), \nonumber \\
\dot{p}_j &=  f_p(V_j,p_j), \nonumber \\
\dot{r}_j &= f_r(V_j,r_j), \nonumber \\
 \quad j &= 1,\dots,N.
\end{align} 
Above, $V_j$ is the transmembrane voltage of neuron $j$, $p_j$ and $r_j$ are associated gating variables, $N$ denotes the total number of neurons,  $u(t)$ is a transmembrane current stimulus common to each neuron, $\sqrt{2 D} \eta_j$ is a white noise process with intensity $D = 1$ associated with neuron $j$, and $C = 1 \mu {\rm F}/{\rm cm}^2$ is a membrane capacitance.  For simplicity, neurons are coupled using all-to-all electrotonic coupling \cite{john95} with strength $\sigma_c$;  other types of neural coupling could also be considered.  Each of the remaining functions from \eqref{neurmod} is described in \cite{rubi04}.  The baseline current $i_b = 5 \mu {\rm A}/{\rm cm}^2$ is chosen so that in the absence of input, coupling, and noise, each oscillator is in a tonically firing regime with a period of $T=8.39$ ms.  In the limit that both $u$, $D$, and $\sigma_c$ are all small in magnitude, Equation \eqref{neurmod} can be well-approximated in a phase reduced form \cite{winf01}, \cite{erme10}, \cite{izhi07}
\begin{align} \label{phasepop}
\dot{\theta}_j &= \omega + Z(\theta_j)\bigg( u(t) + \sqrt{2D} \eta_j + \frac{1}{N} \sum_{i = 1}^N  \sigma_c  (V(\theta_i) - V(\theta_j))  \bigg), \nonumber \\
j &= 1,\dots, N,
\end{align}
where $\theta_j \in [0,2\pi)$ is the phase of oscillator $i$, $\omega = 2\pi/T$ and  $Z(\theta)$ is the phase response curve that characterizes the effect of infinitesimal inputs on the phase.  In the limit as $N \rightarrow \infty$, Equation \eqref{phasepop} can be considered according to a probability density $\rho(\theta,t)$ governed by the Fokker-Planck Equation \cite{gard04}
\begin{equation} \label{fpreduc}
\frac{\partial \rho}{\partial t} =  -\frac{\partial}{\partial \theta} [( \omega + Z(\theta)(u(t) + \sigma_c (\overline{V}-V(\theta)))) \rho(\theta,t)  ] +   \frac{\partial^2}{\partial \theta^2} [ D Z^2(\theta) \rho(\theta,t)],
\end{equation}
with periodic boundary conditions.  Here, $\overline{V} = \int_0^{2\pi} V(\theta) \rho(\theta) d\theta$ is the average voltage.  In previous work \cite{toth22}, \cite{wils20stab}, the above equation was analyzed in the context of developing control strategies to desynchronize a pathologically synchronized population of neural oscillators.    Here, Equation \eqref{fpreduc} will be used in conjunction with the proposed data-driven model identification strategy.  To obtain simulated data for this purpose, the functions $Z(\theta)$ and $V(\theta)$ are computed numerically for the individual neurons from Equation \eqref{neurmod} and Equation \eqref{fpreduc} is subsequently simulated using finite difference approximations for the partial derivatives.  For this model, the state is $\rho_i = \rho(\theta, \Delta t (i-1))$.  The model identification strategy from Section \ref{nlincont} is implemented using 1000 ms of simulated data taking the time step to be $\Delta t = 0.3$ ms.   The input $u(t)$ for training is taken as follows:~random numbers between -1 and 1 are chosen from a uniform distribution with the value held constant over a 2 ms interval; the resulting \added{signal} is smoothed with a spline interpolation and used as the input.

 Two seperate observables are considered for the model \eqref{fpreduc}.  For the first, 
 \begin{equation}  \label{singleobservable}
 g(\rho_i) = \rho(0,\Delta t(i-1)) \in \mathbb{R}^1.
 \end{equation}
To implement the model identification strategy, a delay embedding of size $z = 30$ is used.  No preliminary lifting is considered so that $h(\rho_i) = g(\rho_i)$ as defined in Equation \eqref{hvec}.  As such, $\gamma_{c,i} \in \mathbb{R}^{61}$  as defined in Equation \eqref{gammacont}.  The nonlinear lifting function $f_{c,n}(\gamma_{c,i}) \in \mathbb{R}^{5952}$ is comprised of all possible combinations of polynomial terms taken from $h(x_i), h(x_{i-1}), \dots, h(x_{i-z})$ up to degree 3.  The matrix $\hat{A}$ from Equation \eqref{ahatest} is estimated using a truncated singular value decomposition of rank 40 to approximate the pseudoinverse.  This yields approximations of $A_c \in \mathbb{R}^{61 \times 61}$, $B_c \in \mathbb{R}^{61 \times 1}$, and $C_c \in \mathbb{R}^{61 \times 5952}$ from Equation \eqref{contsys}.  As described in Section \ref{redsec}, the resulting nonlinear equation is projected onto a 20 element POD basis obtained from the eigenvectors of $\Gamma_c \Gamma_c^T$.  Comparisons are also provided  using the Koopman model predictive control strategy from \cite{kord18} which provides a least squares estimate for the update rule $a_i^+ = A a_i + B u_i$ where the lifted state space in this example is taken to be $a_i = \begin{bmatrix} \gamma_{c,i}^T & f_{c,n}(\gamma_{c,i})^T \end{bmatrix}^T$.

\begin{figure}[htb]
\begin{center}
\includegraphics[height=2.2in]{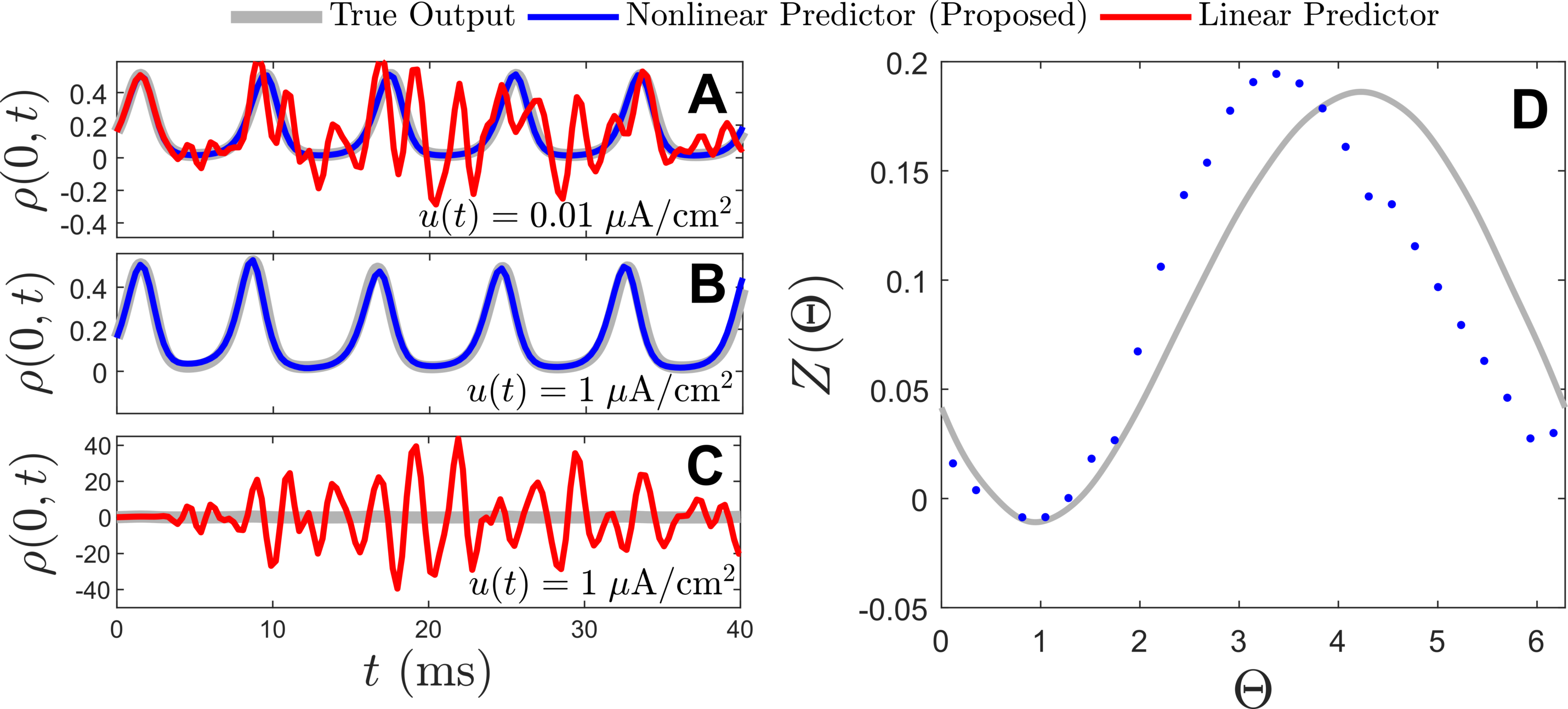}
\end{center}
\caption{Response to pulse inputs for the true model \eqref{fpreduc} and the data-driven inferred models that use either linear or nonlinear predictors. Panel A shows the response to a $0.01 \; \mu {\rm A}/{\rm cm}^2$ pulse starting at 1.5 ms.  The true model output is nearly indistinguishable from the output from the inferred nonlinear model.  By contrast, the inferred linear model displays spurious high frequency oscillations.  Panels B and C show the response to a stronger $1 \; \mu {\rm A}/{\rm cm}^2$ pulse.  For this larger pulse, the nonlinear predictor is still nearly identical to the true model output.   The spurious oscillations from the linear predictor become larger when considering the larger input as shown in panel C.  Panel D compares the population phase response curve from Equation \eqref{phasered} obtained from simulations of the true model (gray line) and the nonlinear inferred model (blue dots).  Note that while pulse inputs were not used for training, the nonlinear inferred model accurately reproduces the response to these pulses.}
\label{neuralpopulation}
\end{figure}

Figure \ref{neuralpopulation} shows simulations using the true model and inferred models in response to pulse inputs.  In panel A, a pulse input of magnitude $0.01 \; \mu {\rm A}/{\rm cm}^2$ lasting 9 ms in duration is applied starting at $t = 1.5$ ms.  Output from the proposed nonlinear predictor is nearly indistinguishable from the output from simulations of the true model \eqref{fpreduc}.  Conversely, the output from the model obtained using the linear predictor is accurate for only the first three milliseconds and ultimately develops spurious high frequency oscillations that render the results inaccurate.  The differences between the linear model and the nonlinear inferred models become more pronounced when using larger inputs; panels B and C show the influence of a magnitude $1 \; \mu {\rm A}/{\rm cm}^2$ pulse with the same timing as the one considered in panel A.  In this case, the nonlinear model still performs well, with outputs that are nearly indistinguishable from the true model outputs.  In response to the $1 \; \mu {\rm A}/{\rm cm}^2$ pulse, the spurious oscillations in the linear model shown in panel C become substantially more pronounced.  

The nonlinear model has stable oscillations when taking $u(t) = 0$ allowing for the further reduction to a phase model of the form \cite{winf01}, \cite{erme10}, \cite{izhi07}
\begin{equation} \label{phasered}
\dot{\Theta} = \Omega + Z(\Theta) u(t).
\end{equation}
Here $\Theta \in [0,2 \pi)$ denotes the phase of the population oscillation, $\Omega$ is the associated natural frequency, and $Z(\Theta)$ is the phase response curve to infinitesimal inputs. Note that capitol Greek letters are used to emphasize that the phase and natural frequency are associated with the population oscillation (as opposed to the phase and natural frequencies of the individual oscillators as considered in Equation \eqref{phasepop}).   Here, $\Theta = 0$ will be defined to occur the moment that $\rho(0,t)$ crosses 0.16 with a positive slope.   $Z(\Theta)$ can be estimated according to the direct method \cite{izhi07}, \cite{neto12}.  This strategy is implemented by applying an pulse input $u(t) = M = 1.5 \mu {\rm A}/{\rm cm}^2$  for a duration $L = 1.5$ milliseconds at a known phase $\Theta_0$ and subsequently inferring the resulting phase shift $\Delta \Theta$. This process is then repeated for different values of $\theta_0$ to provide pointwise estimates of $Z(\Theta_0) \approx \Delta \Theta / M L$ for both the nonlinear inferred model and the true model \eqref{fpreduc}.    Results are shown in panel D of Figure \ref{neuralpopulation} illustrating that the phase response of the true model to inputs is nearly identical to the phase response of the nonlinear inferred model.  Note that it is not possible to provide a similar estimate for the linear inferred model because it does not have a stable periodic orbit.

A second observable is considered for the model \eqref{fpreduc} to illustrate the generality of the proposed model identification strategy.  This second observable is taken to be
\begin{equation} \label{multiobs}
g( \rho_i ) = \begin{bmatrix}  \rho(0,\Delta t(i-1))  \\ \rho(2 \pi/25,\Delta t(i-1)) \\ \rho(4 \pi / 25,\Delta t(i-1)) \\ \vdots \\ \rho(48 \pi/25,\Delta t(i-1))  \end{bmatrix} \in \mathbb{R}^{25},
\end{equation}
i.e.,~the observable is comprised of 25 measurements of $\rho(\theta,t)$ equally spaced in $\theta$.  With this alternative observable,  the model identification strategy is implemented using a delay embedding of size $z = 20$.  Data for model identification is taken for 1000 ms of simulated data with a time step of $\Delta t = 0.3$ ms.   The input $u(t)$ used for training is chosen as follows:~random numbers between -0.25 and 0.25 are chosen from a uniform distribution with the value held constant over a 2 ms interval.  The resulting input is smoothed with a spline interpolation and applied for simulations of the true model \eqref{fpreduc}.  Once again, no preliminary lifting is considered  so that $h(\rho_i) = g(\rho_i)$.  Consequently $\gamma_{c,i} \in \mathbb{R}^{545}$.  The nonlinear lifting function $f_{c,n}(\gamma_{c,i}) \in \mathbb{R}^{3250}$  is comprised of all possible combinations of  polynomial terms taken from $h(\rho_i)$ up to degree 4.  The matrix $\hat{A}$ is estimated according to Equation \eqref{ahatest} using a truncated singular value decomposition of rank 40 to approximate the pseudoinverse.  This information is used to determine $A_c \in \mathbb{R}^{545 \times 545}$, $B_c \in \mathbb{R}^{545 \times 1}$, and $C_c \in \mathbb{R}^{545 \times 3250}$ from Equation \eqref{contsys}.  As described in Section \ref{redsec}, a 25-dimensional model is obtained by projecting the inferred model equations onto a POD basis obtained from the eigenvectors of $\Gamma_c \Gamma_c^T$.  Once again, comparisons are provided when using the Koopman model predictive control strategy \cite{kord18} that obtains a least squares estimate for the update rule $a_i^+ = A a_i + B u_i$ where the lifted state space in this example is taken to be $a_i = \begin{bmatrix} \gamma_{c,i}^T & f_{c,n}(\gamma_{c,i})^T \end{bmatrix}^T$.  

\begin{figure}[htb]
\begin{center}
\includegraphics[height=3in]{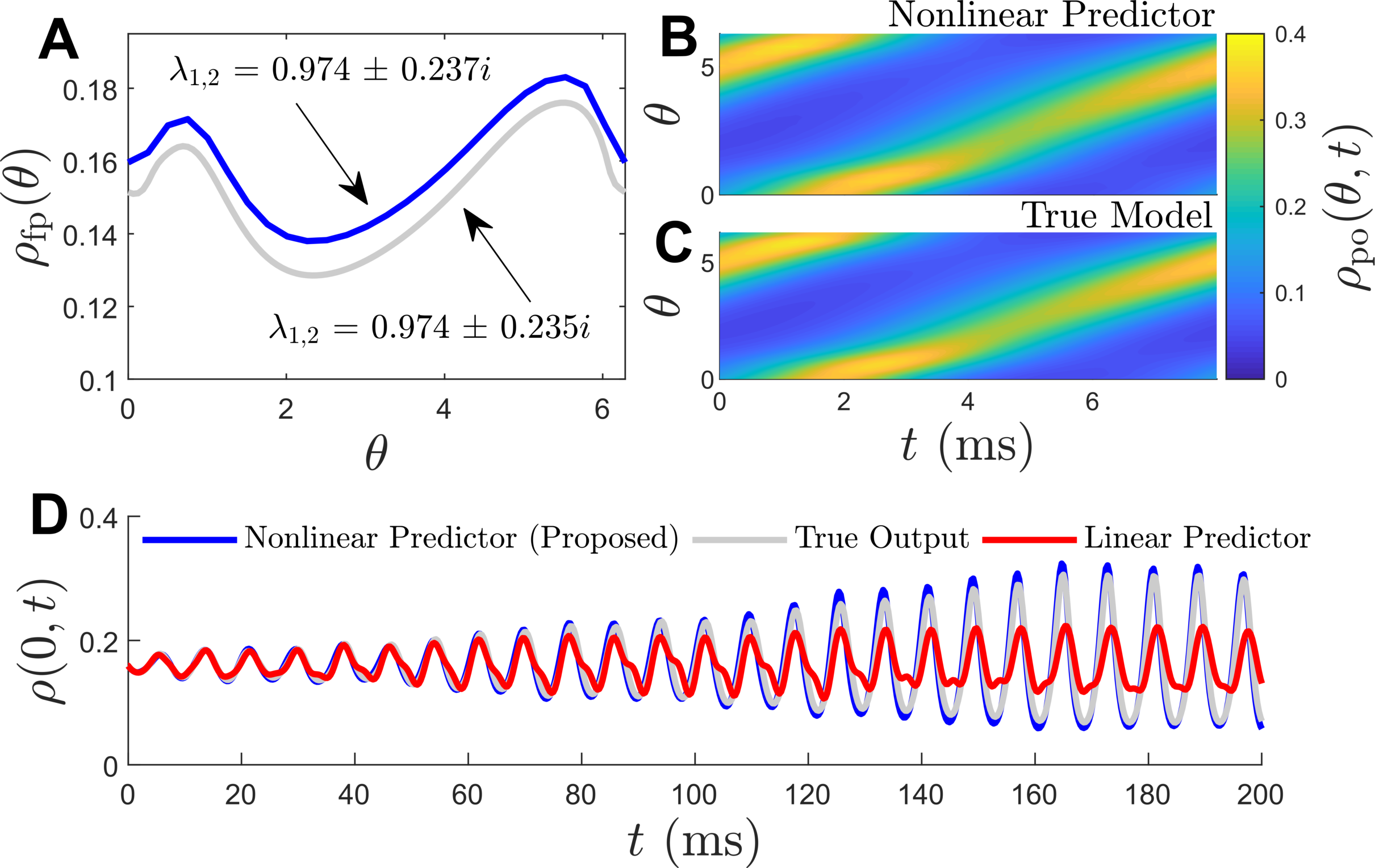}
\end{center}
\caption{Dynamical features of the model obtained using the nonlinear predictor with observables from \eqref{multiobs}.  Panel A shows an unstable fixed point, $\rho_{\rm fp}(\theta)$, for the inferred nonlinear model (blue line) which is nearly identical to the unstable fixed point of the full model (gray line).  Unstable eigenvalues associated with these solutions are nearly identical.  Initial conditions near the unstable fixed points of each model are in the basin of attraction of a stable periodic orbit.  Panels B and C show that the periodic orbits, $\rho_{\rm po}(\theta,t)$, obtained for each model are nearly identical.  Panel D highlights a trace of $\rho(0,t)$ for an initial condition near the stable fixed point.  Output from the nonlinear predictor matches the output from the true model.  Output from the linear model matches well in the early parts of this simulation, but does not capture the full increase in the amplitude of the oscillation.}
\label{neuralpopfulldata}
\end{figure}

In addition to accurately predicting the output of the true model in response to input, the inferred nonlinear model accurately characterizes fixed points and periodic orbits as shown in Figure \ref{neuralpopfulldata}.  For instance, panel A illustrates an unstable fixed point, $\rho_{\rm fp}(\theta)$, that exists both in the full model \eqref{fpreduc} and the nonlinear inferred model of the form \eqref{contsys}  when taking $u(t) = 0$.  As indicated in the figure, both the profile of the fixed point solution and the the associated unstable, complex-conjugate discrete time eigenvalues (obtained from the linearization about the fixed point) are nearly identical.  This unstable fixed point emerges in the true model as a result of a Hopf bifurcation where the coupling strength is the bifurcation parameter.     Note that the model obtained when using the observable \eqref{singleobservable} also has a fixed point with $\rho(0,t) = 0.159$ with associated unstable discrete time eigenvalues $\lambda_{1,2} = .9740 \pm 0.239 i$.    Despite using different data sets, the models inferred from the from the observables \eqref{singleobservable} and \eqref{multiobs} yield comparable estimates for the location and stability of this fixed point.    Of course, for the model that uses the single observable \eqref{singleobservable}, it is not possible to reconstruct the full probability density since this information is unavailable.  Initial conditions near this unstable fixed point eventually settle to a stable periodic orbit; the periodic orbits, $\rho_{\rm po}(\theta,t)$, obtained from the nonlinear predictor and the true model represented by the colormaps in panels B and C, respectively, are nearly indestinguishable.  Panel D illustrates the time course of $p(0,t)$ for an initial condition near the unstable fixed point when taking $u(t) = 0$.  The transition from the stable fixed point to the unstable periodic orbit is well captured by the proposed nonlinear predictor.  The model obtained using the Koopman model predictive control strategy \cite{kord18} does not accurately reflect this transition.

\FloatBarrier

\subsection{One Dimensional Burgers' Equation}

The Burgers' equation is often used as a test bed for Koopman-based model identification and analysis strategies \cite{page18}, \cite{arba18b}, \cite{peit19}   because it has a convective nonlinearity that is similar to that of the Navier-Stokes equations.  Here a 1-dimensional version of the Burgers' equation is considered
\begin{equation} \label{burgeq}
\frac{\partial w}{\partial t} = \frac{1}{{\rm Re}} \frac{\partial^2 w}{\partial x^2} - w \frac{\partial w}{\partial x}.
\end{equation}
Here $w(x,t)$ gives the state on the domain $x \in[0,1]$ and ${\rm Re} = 50$ is a constant that is analogous to the Reynolds number from the Navier-Stokes equations.  In this example, Dirichlet boundary conditions $w_L(t)$ and $w_R(t)$ are considered for the boundary at $x = 0$ and $x = 1$, respectively.  These boundary conditions are also taken to be the inputs, i.e.,~$u(t) = \begin{bmatrix} w_L(t) & w_R(t)  \end{bmatrix}^T$. For this model, the state is $w_i = w(x, \Delta t(i-1))$.  The model identification strategy from Section \ref{nlincont} is implemented taking 2000 time units of simulated data with $\Delta t = 0.1$ ms.  The input $u(t)$ used for training is chosen as follows:~for both $w_L(t)$ and $w_R(t)$, random numbers between -0.5 and 0.5 are chosen from a uniform distribution with the value held constant over a 20 time unit interval.  These signals are smoothed with a spline interpolation and the resulting inputs are used in training simulations.

The observable for the model \eqref{burgeq} is taken to be
\begin{equation} \label{burgoutput}
g(w_i) = \begin{bmatrix} w(0, \Delta t(i-1)) \\ w(0.05, \Delta t(i-1))  \\ w(0.10, \Delta t(i-1)) \\ \vdots  \\ w(0.95, \Delta t(i-1)) \end{bmatrix} \in \mathbb{R}^{20},
\end{equation}
i.e.,~the observable is comprised of 20 measurements of $w(x,t)$ equally spaced in $x$.  Following the definitions given in Section \ref{koopbackground} and \ref{koopnonlin}, a delay embedding of size $z = 30$ is used.  No preliminary lifting is considered so that $h(w_i) = g(w_i)$ as defined in Equation \eqref{hvec}.  Here $\gamma_{c,i} \in \mathbb{R}^{680}$ as defined in Equation \eqref{gammacont}.  The nonlinear lifting function $f_{c,n}(\gamma_{c,i}) \in \mathbb{R}^{1750}$ is comprised of all possible combinations of polynomial terms taken from $h(w_i)$ up to degree 3.  The matrix $\hat{A}$ from Equation \eqref{ahatest} is estimated using a truncated singular value decomposition of rank 80 to approximate the pseudoinverse.  This yields approximations of $A_c \in \mathbb{R}^{680 \times 680}$, $B_c \in \mathbb{R}^{680 \times 2}$, and $C_c \in \mathbb{R}^{680 \times 1750}$ from Equation \eqref{contsys}.  As described in Section \ref{redsec}, the resulting nonlinear equation is projected onto a lower dimensional POD basis obtained from the eigenvectors of $\Gamma_c \Gamma_c^T$.  Comparisons are also provided  using the Koopman model predictive control strategy from \cite{kord18} which gives an estimate for the update rule $a_i^+ = A a_i + B u_i$ where the lifted state space in this example is taken to be $a_i = \begin{bmatrix} \gamma_{c,i}^T & f_{c,n}(\gamma_{c,i})^T \end{bmatrix}^T$.  This least squares fitting is implemented according to Equation \eqref{dmdcfit} using a truncated singular value decomposition retaining different numbers of singular values as described in the results below.

Comparisons are provided between simulations of the true model \eqref{burgeq} and both the linear and nonlinear inferred models in response to different inputs.  Panels A and B show the $L^2$ error, $L^2 = \int_0^1 (w_{\rm true}(t,x)- w_{\rm inferred}(t,x))^2dx$ between the true model solutions and the inferred solutions when using the linear predictor (red lines) and nonlinear predictor (blue lines) for different inputs.  In panel A, $w_L(t)$ and $w_R(t)$ are chosen similarly to the training data, i.e.,~obtained by choosing random numbers between -0.5 and 0.5 from a uniform distribution, holding the value constant over a 20 time unit interval, and smoothing the resulting input with a spline interpolation.  Note that the inputs are not identical those used for training because the random numbers are realized differently.  Results from panel B consider a similar input, except that the values that comprise $w_L$ and $w_R$ are held constant for only 7 time units before smoothing.  Effectively, this yields inputs with higher frequency content than those used for training.  In each case, the $L^2$ error is approximately three orders of magnitude lower when using the nonlinear predictor as compared to the linear predictor.  Panels C and D show representative traces of $w(x,t)$ corresponding to simulations in panels A and B, respectively.  The true model output and the output from the nonlinear model are nearly identical in both cases while the linear predictor does not accurately capture the model output.  Panel E shows the average $L^2$ error for inferred models of different order when applying the lower frequency input.  For the linear model, order is determined by the number of singular values retained in the truncated singular value decomposition used to approximate the pseudoinverse in Equation \eqref{dmdcfit}.  For the nonlinear model, order is governed by the dimension of the POD basis used for projection as described in Section \ref{redsec}.  The linear model performs slightly better as the order increases.  For moderate orders between 10 and 100, the inferred linear model is generally unstable, i.e.,~the inferred matrix $A$ from Equation \eqref{statecontrol} has unstable eigenvalues.  These dots are omitted from panel E because the error grows unbounded as time approaches infinity.  The inferred nonlinear models do not suffer from the same stability issues as the linear models.   Accuracy of the nonlinear model stops improving once the model order reaches 50, at which point the output is nearly indistinguishable from the true model output.

\begin{figure}[htb]
\begin{center}
\includegraphics[height=2.6in]{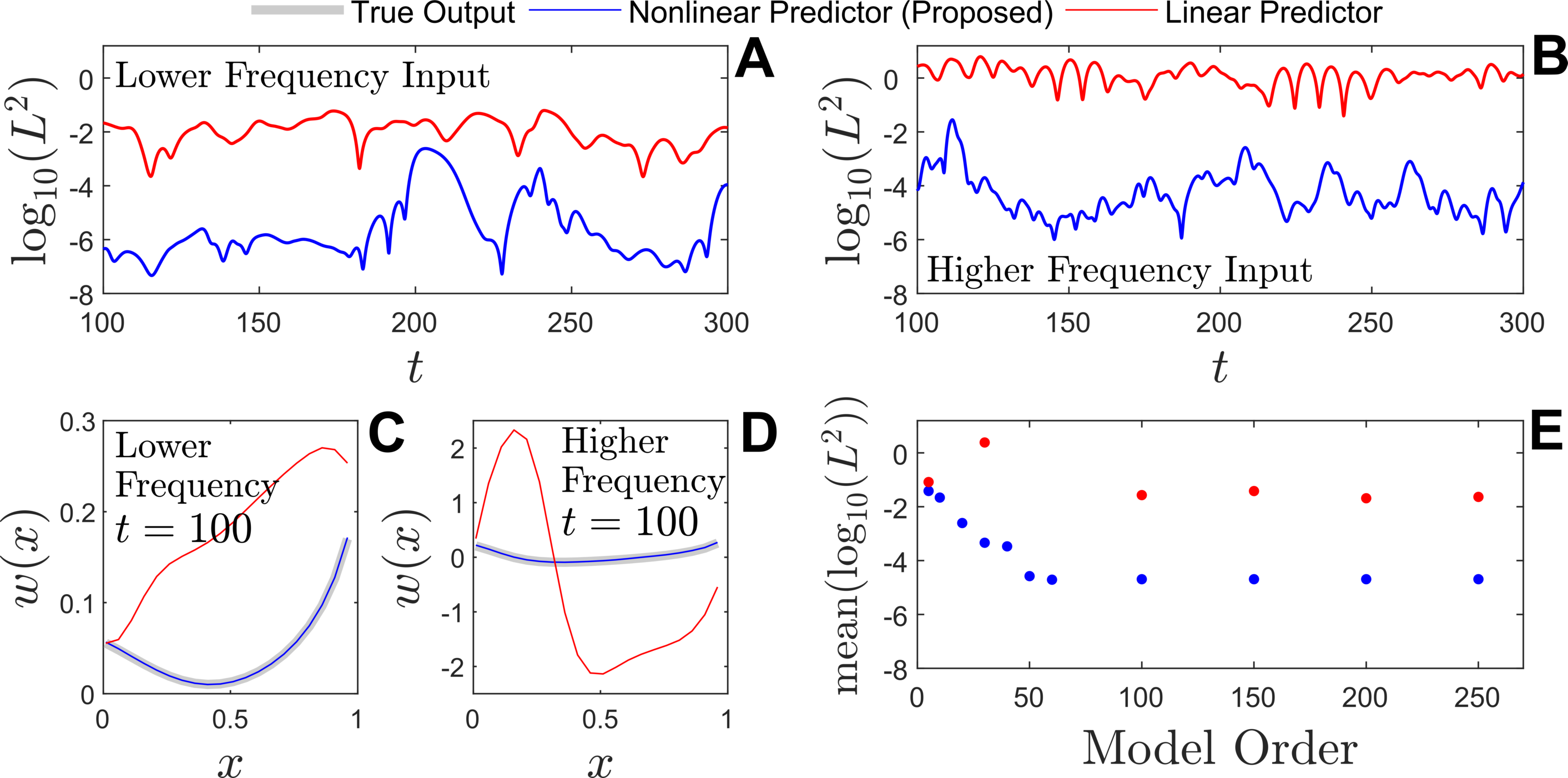}
\end{center}
\caption{Accuracy of the linear and nonlinear models inferred from data in response to inputs with different frequency content.  In panel A, inputs are similar to those used for training are considered.  In panel B, inputs with higher frequency content are considered.  These panels show representative traces of the $L^2$ error over a 200 time unit window of simulation. In each case, the proposed nonlinear predictor yields results that are approximately 3 orders of magnitude better than the linear predictor.  Panels C and D give representative traces of the outputs from the simulations from panels A and B, respectively.  In each case, the nonlinear predictor yields outputs that are nearly identical to the true model output while the linear predictor yields outputs that are substantially less accurate.  Panel E shows the influence of the model order on the accuracy of the inferred linear and nonlinear models.  For models with order between 10 and 100, most inferred linear models are unstable with errors that grow unbounded in time.  As such, there are fewer data points for the linear predictor in panel E.}
\label{burgersresults}
\end{figure}

\subsection{Schlierin Images of Supersonic Flow Past a Cylinder} \label{schsec}
\begin{figure}[htb]
\begin{center}
\includegraphics[height=2.4 in]{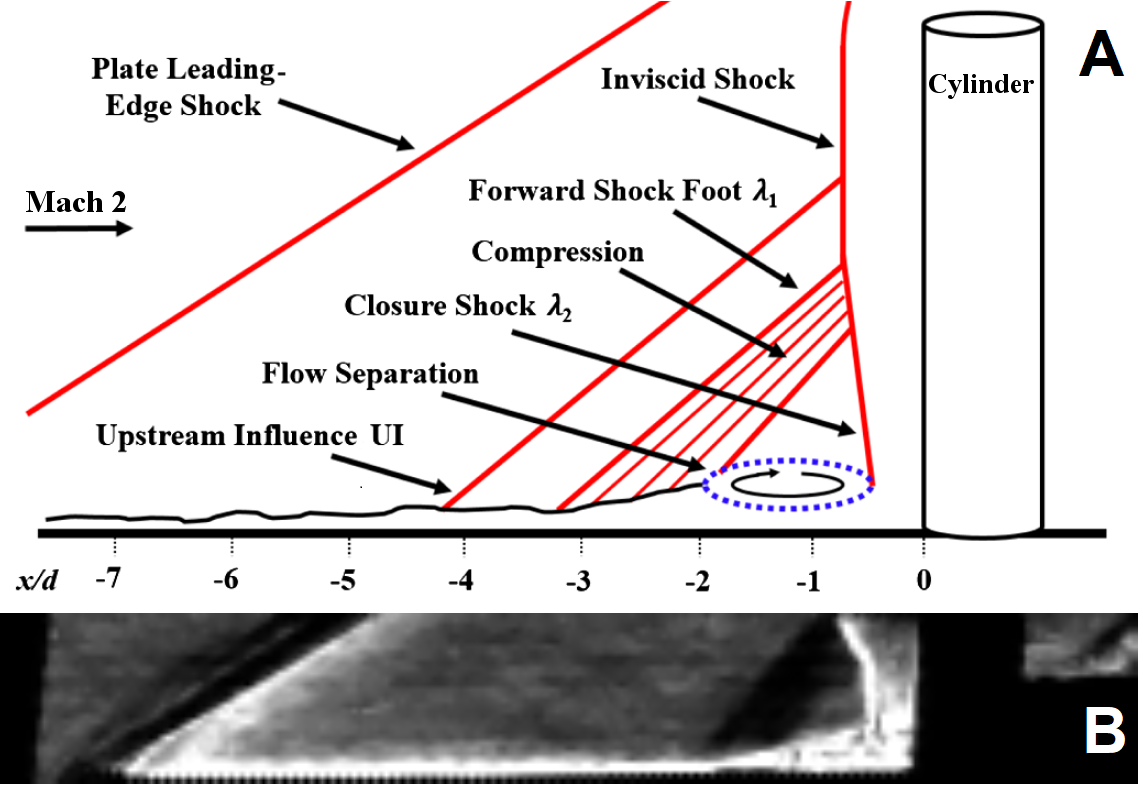}
\end{center}
\caption{Panel A shows a schematic depicting the flow geometry used to study shock-wave/transitional boundary-layer interaction.  Mach 2 flow enters from the left and interacts with the cylinder mounted to a flat plate.  Temporal oscillation in the location of the forward shock foot is of particular interest here.  Panel B shows a characteristic schlieren image taken from this data set.  \added{Adapted from \cite{wils20acc}.}}
\label{lambda}
\end{figure}

Finally, experimental schlieren image data of cylinder-generated shock-wave/transitional boundary-layer interaction is analyzed using the proposed nonlinear model identification strategy.  Here, the schlieren images of the Mach 2 flow past a standing cylinder are taken at 100 kHz.  Details of the experimental setup and data collection are provided in \cite{wils20acc}.   Salient flow features are illustrated in panel A of Figure \ref{lambda} with the flow going from left to right.    Panel B shows a characteristic schlieren image taken from this data set.  A flat plate is visible on the bottom edge and the cylinder is visible near the right edge of the image.  Differences in pixel intensities roughly correspond to differences in fluid density gradients.  Of particular interest in this data is the location of the forward lambda shock foot.  Previous studies identified a characteristic oscillation frequency of the forward shock foot of approximately 5 kHz \cite{comb18}, \cite{comb19}, \cite{wils20acc} by analyzing power spectral densities and by using linear data analysis techniques such as spectral POD and other techniques related to dynamic mode decomposition.  

The data set consists of 25,000 snapshots with each image containing of 5,472 pixels.  In order to make the data set more computationally tractable, every other snapshot is removed so that the data is effectively sampled at 50 kHz.   Half of the remaining snapshots are used for training and the other half are used to validate the resulting data-driven model.   To further compress the data, snapshots from the training set are represented using a 5 mode POD basis which captures 0.48 of the total energy as determined by taking the sum of the largest 5 eigenvalues of the covariance matrix and dividing by the total sum of the eigenvalues.   The observable is taken to be
\begin{equation} \label{podobs}
g_i = \begin{bmatrix}  \omega_{1,i} & \dots & \omega_{5,i} \end{bmatrix}^T \in \mathbb{R}^5,
\end{equation}
where $\omega_{k,i}$ is the amplitude of the $k^{\rm th}$ POD mode on the $i^{\rm th}$ snapshot.  A representation in the space of the schlieren images can be obtained by taking a linear combination of the POD modes with weights $\omega_{1,i}, \dots, \omega_{5,i}$.  The strategy from Section \ref{autsys} is implemented on the autonomous data set in order to infer a nonlinear model.  Once again, no initial lifting is considered so that $h_i = g_i$. A delay embedding of size $z = 25$ is used.  Consequently, $\gamma_{i} \in \mathbb{R}^{130}$ as defined in Equation \eqref{liftstate}.  The nonlinear lifting function $f_n(\gamma_i) \in \mathbb{R}^{120}$ is comprised of all possible combinations of polynomial terms that comprise $h_i$ up to degree 4.  The matrix $\hat{A}$ is estimated according to Equation \eqref{minsol} yielding approximations of $A_n \in \mathbb{R}^{130\times 130}$ and $C_n \in \mathbb{R}^{130 \times 120}$ in the nonlinear estimator from Equation \eqref{nonlinpredict}.  As described in Section \ref{redsec}, the resulting nonlinear equation is projected onto a low rank basis comprised of the 4 most important POD modes obtained from the eigenvectors of $\Gamma \Gamma^T$.

\begin{figure}[htb]
\begin{center}
\includegraphics[height=5.0 in]{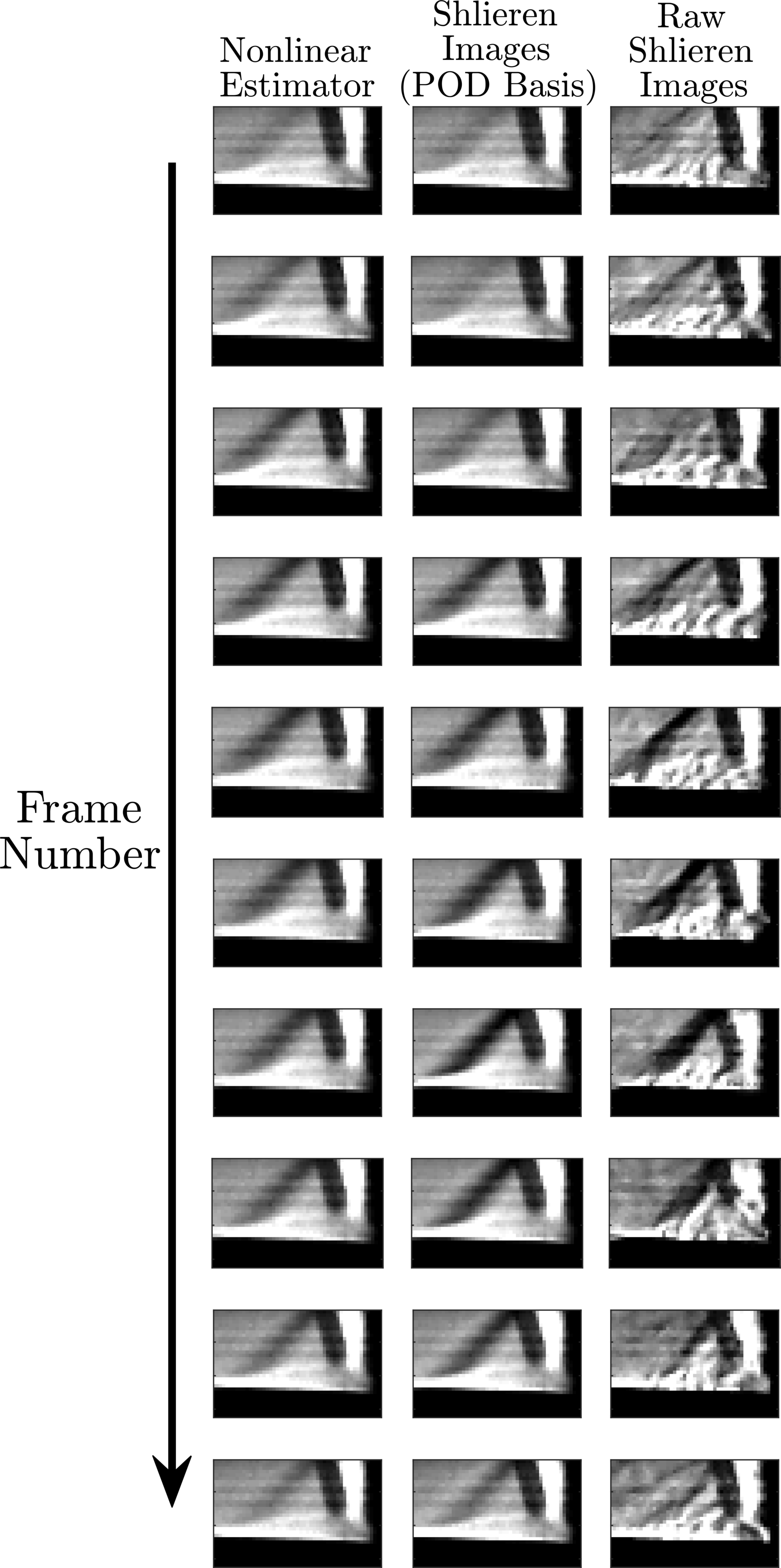}
\end{center}
\caption{The inferred nonlinear model of the form \eqref{nonlinpredict} has a stable periodic orbit with frequency 4.70 kHz shown in sequential frames in the left columns.  In the fourth frame from the top, the forward shock foot is farthest from the cylinder.  The forward shock foot is closest to the cylinder in the eighth and ninth frames.  The right columns show sequential schlieren images from the comparison data set (i.e.,~that was not used during training).  The middle panel shows these images projected onto the same 5 POD mode basis used to implement the nonlinear model identification strategy.  While the left and middle columns display slight variations in the pixel intensities, the periodic orbit identified by the inferred nonlinear model is nearly indistinguishable from the behavior observed in the experimental data.}
\label{shockwaveresults}
\end{figure}

The inferred nonlinear model of the form \eqref{nonlinpredict} has a stable periodic orbit with frequency of 4.70 kHz that which agrees well with the oscillation frequencies identified in prior studies \cite{comb18}, \cite{comb19} \cite{wils20acc}.  This periodic orbit is identified by taking an initial condition obtained from the comparison data set (i.e.,~the portion of the data set not used for training) and iterating the model \eqref{nonlinpredict} until the initial transients decay.     In contrast to the results of this study, the linear techniques applied in these previous studies did not explicitly identify a periodic orbit, but rather, identified characteristic oscillation frequencies observed with the data set.  Figure \ref{shockwaveresults} shows the periodic orbit obtained from the inferred nonlinear model (left column) as well as raw data from the comparison data set in the right columns.  The ten sequential frames are centered at the lambda shock to emphasize the oscillation in the location of the forward shock foot and correspond to approximately one period of oscillation.  The cylinder and the flat plate appear in the right and bottom edges each frame, respectively.  The middle column of Figure \ref{shockwaveresults} shows the raw data projected onto the 5 POD mode basis used for the nonlinear model identification strategy; the left and middle frames are qualitatively similar to each other both in terms of the location of the forward lambda shock and in terms of the pixel intensities. Note that the higher frequency flow features, for instance that appear in the flow separation region between the forward shock foot, $\lambda_1$, and the closure shock, $\lambda_2$, are not accurately resolved in the reduced order model; these features are filtered out in the initial projection of the data onto the 5 mode POD basis.

\begin{figure}[htb]
\begin{center}
\includegraphics[height=2.0 in]{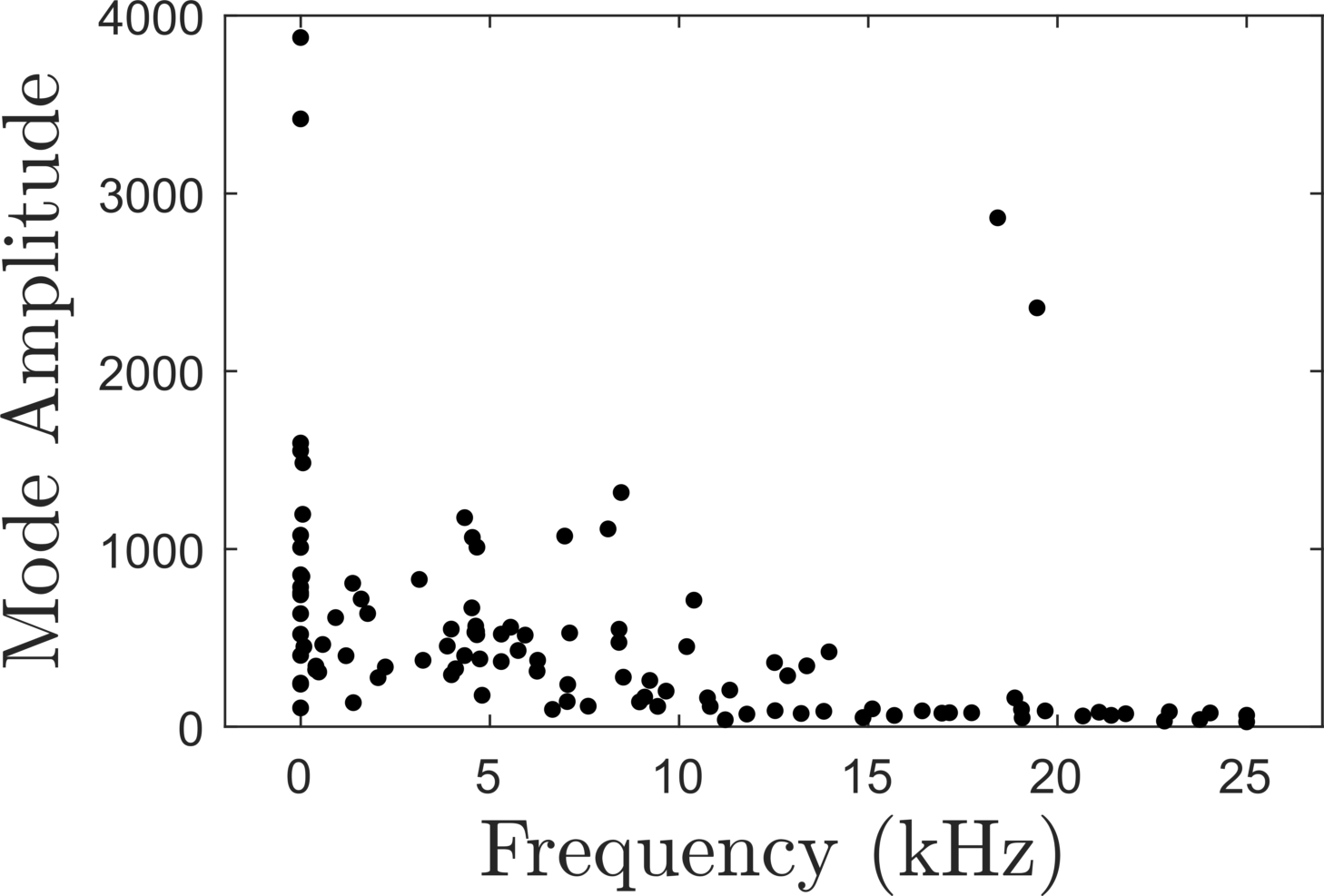}
\end{center}
\caption{ Extended DMD is applied to the schlieren image data.  The resulting eigenmodes are shown according to their frequency content and their relative importance as gauged by the average amplitude of the eigenmode observed in the snapshot data.  Previous analysis  of this data in \cite{comb18} identified peaks in the power spectrum at approximately 4.7 kHz, which is consistent with the stable 4.7 kHz periodic orbit that results when using the proposed nonlinear model identification technique.  By contrast, the extended DMD algorithm identifies two dominant eigenmodes at approximately 20 kHz; the eigenmodes with frequency content near 4.7 kHz do not stand out relative to the other eigenmodes.}
\label{shockwaveedmd}
\end{figure}

For comparison, the extended DMD approach \cite{will15} is also implemented on the schlieren image data set as described in Section \ref{koopest} to provide a linear least squares estimate for the update rule $a_i^+ = A a_i$.  Here, the lifted state space is taken to be $a_i = \begin{bmatrix} \gamma_i^T & f_n(\gamma_i)^T \; \end{bmatrix}^T$.  The least squares fitting uses a truncated singular value decomposition of rank 200 -- keeping more singular values results in an unstable linear system.  This extended DMD approach is often used to provide an approximation for Koopman eigenmodes.  However, as noted in \cite{comb19} it is often difficult to identify which Koopman eigenmodes are most important in a given data set.  Indeed, Figure \ref{shockwaveedmd} shows a plot of frequency of the eigenmodes versus the average amplitude of the eigenmode from the snapshot data.  Here, the frequency associated with a given eigenmode is equal to \added{the imaginary component of} $\log(\lambda_A)/2 \pi \Delta t$ where $\lambda_A$ is an eigenvector of the inferred $A$ matrix and $\Delta t$ is the time between successive snapshots.  The associated amplitude at frame $i$ is given by  $w_A^T a_i$ where $w_A^T$ is the left eigenvector for eigenvalue $\lambda_A$.  Analysis of the frequency content of these eigenmodes would identify two dominant eigenmodes with a frequency near 20 kHz.  However, power spectrum analysis of the same schlieren image data performed in \cite{comb18} identifies a peak in power at approximately 4.7 kHz and relatively little power in the 20 kHz region.  By contrast, the proposed nonlinear model identification technique identifies a 4.70 kHz stable periodic orbit embedded in the data which is consistent with peaks in the power spectrum of the imaging data identified in \cite{comb18}.

\FloatBarrier

\section{Conclusion} \label{concsec}
Koopman analysis and associated data-driven model identification algorithms are invaluable tools in the study of nonlinear dynamical systems.  In data-driven modeling applications, the vast majority of Koopman-based strategies consider finite-dimensional, linear estimators for the action of the Koopman operator on observables.  This work proposes a general strategy to obtain a nonlinear estimator for the Koopman operator.  In the examples considered in this work, the proposed strategy yields nonlinear models that are substantially more accurate on longer timescales than comparable models that consider linear estimators for the Koopman operator.  It should be noted that the examples considered in this work did not explicitly consider optimization of the observables used for fitting, for instance, which might result in a low-dimensional Koopman invariant subspace \cite{take17}, \cite{brun16}, \cite{kord20}.  As such one cannot rule out the possibility that a linear estimator could have provided a more accurate representation for the Koopman operator with a more careful choice of the observables for the examples considered in this work.  

Despite the promising results presented in the applications considered in this paper, the proposed nonlinear model identification strategy certainly comes with its share of drawbacks as discussed below.  Foremost, the proposed model identification strategy yields a nonlinear model, thereby precluding the direct use of a wide variety of linear model analysis techniques and control algorithms that are available for linear models obtained using DMD \cite{kutz16}, Extended DMD \cite{will15}, and Koopman model predictive control \cite{kord18}.  In applications where the dynamics can be well-approximated by a low-dimensional linear operator, linear estimators would certainly be preferable to nonlinear estimators.  The proposed algorithm shares similarities with the extended DMD algorithm in that it uses a dictionary of functions of the observables to lift to a higher dimensional space.  Strategies for determining an optimal choice for the lifting functions were not considered here.  Practically, using polynomial combinations of the observables and radial basis functions worked well in the applications presented in this work.  The proposed model identification technique shares similarities with other approaches designed to determine the model equations directly from data.  In contrast to sparse model identification algorithms \cite{mang19}, \cite{brun16b}, \cite{pant19}, \cite{rudy17}, the proposed strategy instead identifies low rank approximations for the nonlinear dynamics.  In contrast to these sparse model identification strategies, the proposed methodology yields a nonlinear  model that is generally not interpretable in the sense that the terms of the learned equation correspond to the physics of the underlying models.  Nonetheless, the proposed strategy does not require the use of neural networks and still provides a low-rank approximation for the dynamics.  

Due to its similarity to existing Koopman-based model identification strategies, the proposed framework offers some interesting opportunities for extension.  For instance, this strategy could be used alongside linear Koopman model predictive control strategies described in \cite{arba18b}, \cite{kord18}, using the linear estimator to approximate an optimal control input over a finite time horizon and subsequently using the nonlinear predictor to maintain accurate information about the state of the system.  Such an approach would be particularly useful in situations where real-time information about the observables is not continuously available.   Alternatively, the resulting nonlinear models can be used to obtain additional information about stable attractors in the inferred systems.  This point was illustrated in Section \ref{neurpopsec} where applying pulse inputs to the nonlinear model allowed for accurate estimation of the phase response curve for the limit cycle that emerges in the true model equations.  The diversity of examples considered in this work using both computational and experimental data suggest that the proposed framework could be a versatile tool to aid in the identification of nonlinear dynamical systems, especially in applications where linear predictors alone are not sufficient.  

\section*{Acknowledgments}
This material is based upon work supported by the National Science Foundation Grant CMMI-2140527.  Thank you to Phil Kreth from University of Tennessee Space Institute for providing the schlieren imaging dataset considered in Section \ref{schsec}.

\end{document}